\documentclass[10pt,twoside]{article}
\usepackage{amsfonts,amsmath,amssymb}
\usepackage{latexsym}

\newtheorem{theorem}{Theorem}[section]

\newtheorem{lemma}{Lemma}[section]
\newtheorem{proposition}{Proposition}[section]

\newtheorem{example}{Example}[section]
\newtheorem{corollary}{Corollary}[section]
\def\proof{\mbox {\it Proof.~}}
\def\a{\alpha}
\def\sh{\sharp}
\def\e{\epsilon}
\def\D{\Delta}
\def\n{\nabla}
\def\U{\Upsilon}
\def\o{\theta}
\def\w{\omega}

\def\d{\delta}
\def\s{\sigma}
\def\m{\mathcal}

\def\R{\mathbb R}
\def\mE{\mathcal E}
\def\Bodeux{B_{0,G_2}(M,g)}     \def\Boun{B_{0,G_1}(M,g)}
\def\Bog{B_0(M,g)}
\def\BoG{B_{0,G}(M,g)}
\def\H{H_1^2(M)}               
       
\def\Kdeux{\frac{K_{n-k}}{A_2^\frac{2}{n-k}}}
\def\Bodeux{B_{0,G_2}(M,g)}
\def\mHun{(\mathcal H_1)}    \def\mHdeux{(\mathcal H_2)}
\def\t{\tilde}

\makeatletter\def\theequation{\arabic{section}.\arabic{equation}}\makeatother

\begin{document}
\fontfamily{ptm}\fontseries{sb}\selectfont
\title{
\begin{flushleft}
{\small{\it Advanced Nonlinear Studies} {\bf 8} (2008), 303--326}
\end{flushleft}
\vspace{0.5in} {\bf\Large  Multiplicity for Critical and
Overcritical Equations}}
\author{{\bf\large Dellinger \,Marie}\vspace{1mm}\\
{\it\small Equipe G{\'e}om{\'e}trie et Dynamique}\\
{\it\small Institut de  Math{\'e}matique, Universit{\'e} Paris VI}\\
{\it\small 175 rue Chevaleret, 75013 Paris, France}\\
{\it\small e-mail: dellinge@math.jussieu.fr}\vspace{1mm}}
\date{\small{Received 23 September 2008} \\
{\it\small Communicated by Jean Mawhin}}

\maketitle
\begin{center}
{\bf\small Abstract}

\vspace{3mm}
\hspace{.05in}\parbox{4.5in}
{{\small
Let  $(M,g)$ be  a compact Riemannian  $n-$manifold, $\ n\geq 3.$
 We prove the existence of multiple solutions for equations like
$$
\D u + \a u = f u^{p}, \qquad u>0
$$
where $\a\in \R^{+ *}, f\in C^\infty(M)$ is positive, and the
exponent $p$ takes critical and overcritical values.
General results are obtained and specific examples are discussed,
 like  $S^n, S^1(t)\times S^{n-1},$ and
$S^1(a)\times S^2(b) \times S^{n-3}.$}}
\end{center}

\noindent
{\it \footnotesize 2000 Mathematics Subject Classification}. {\scriptsize 58J05 (35B33)}.\\
{\it \footnotesize Key words}. {\scriptsize Nonlinear elliptic
  equations, Riemannian manifold, Sobolev inequality, Yamabe problem.}

\section{\bf Introduction}
\def\theequation{1.\arabic{equation}}\makeatother
\setcounter{equation}{0}

Let $(M,g)$  be a compact Riemannian  manifold of dimension $n\geq 3.$
Our  paper is concerned with the question of the existence of
multiple smooth  solutions  for the  equation
\begin{equation}\tag{$E_p$}
\D u + \a u = f u^p, \qquad u>0
\end{equation}
where $\D = - \mbox{div }  ( \nabla \  )   $ is the $g-$Laplacian,
 $\a\in \R^{+ *}, f\in C^\infty (M)$ is  positive and $p\geq
\frac{n+2}{n-2}.$ We say that the equation $(E_p)$ is critical when
$p=\frac{n+2}{n-2}$ and overcritical when $p>\frac{n+2}{n-2}.$ Indeed, the
exponent $\frac{n+2}{n-2} $ is the classical critical Sobolev growth
exponent. It appears in particular in the equation one has to solve
in the prescribed scalar curvature problem :
\begin{equation}\label{EScal}
\D u + \frac{(n-2)\  S_g}{4(n-1)} u = f u^\frac{n+2}{n-2},  \qquad u>0
\end{equation}
where $S_g $ is the scalar curvature of $g.$
More precisely, if for $f\in C^\infty(M)$ there exists $u\in
C^\infty(M)$ a positive solution of (\ref{EScal}), then $f$ is the
scalar curvature of the $g$-conformal metric $ u^\frac{4}{n-2} g.$
We are here    interested in two particular cases of  equation (\ref{EScal}).
On the standard sphere $(S^n,h_n),$ this problem  is referred to as
the Nirenberg problem. Its resolution is equivalent
to the resolution of  (\ref{EScal}) with $S_{h_n} = n(n-1).$
For references on the Nirenberg problem,  see Hebey \cite{Hnir},
Kazdan-Warner \cite{KW} and Li \cite{L}.
 There is also the intensively studied  Yamabe problem, which consists
 in  the search for  conformal metrics  with constant scalar
 curvature.
It corresponds to the resolution of (\ref{EScal}) with
$f=1.$ The Yamabe problem is completly solved. \\

  Concerning multiplicity and uniqueness of positive solutions for
such equations,  we refer  to Aubin \cite{A, Abook},
Bidaut-V{\'e}ron and V{\'e}ron \cite{BVV}, Esposito \cite{E},
Hebey-Vaugon \cite{HVmult}, Obata \cite{O}, Pollack \cite{P}, Schoen
\cite{Schoen1} and  \cite{Schoen2}. In particular, note that the
Yamabe equation possesses a unique
 solution if there exists $\t g \in  [g]$  such that $S_{\t g}\leq 0$
or if there exists an Einstein metric $\t g
\in [g],$ where $[g]$ stands for the conformal class of $g.$
 We are here especially interested on results of  Hebey-Vaugon
 \cite{HVmult}
(see also  Schoen \cite{Schoen1}).
In their work, the manifold is assumed to have big enough  isometry groups  and
 solutions are required to be invariant under the action of subgroups.
Besides, all groups are finite which implies that the quotient space
of all orbits can be equiped with a structure of manifold.
In our results, this  condition is not required. This is made possible
thanks to the recent advances of Hebey-Vaugon \cite{HVsobosym} and Faget
\cite{zoe03,zoe04} concerning the influence of isometry groups
 on Sobolev spaces and Sobolev inequalities. \\

 Given  $G$  an isometry group, $\a\in \R^{+*},$ and $f\in
C^\infty(M)$
 positive and $G$-invariant, we consider
$G$-invariant solutions of the equation
\begin{equation}\tag{$E_{\a f }^k$}
\D u + \a u = f u^\frac{n+2-k}{n-2-k}, \qquad u>0,
\end{equation}
where $k\geq 0$ is the minimum dimension of the $G$-orbits.
The energy of a  solution $u$  of $(E_{\a f}^k)  $ is defined by
\begin{equation}\label{defnrj}
\mE(u)= \int_M f u^\frac{2(n-k)}{n-2-k} \ dv_g.
\end{equation}
We obtain multiplicity of energies for solutions of $(E_{\a f}^k)  $ where
 each solution is invariant by the action of an isometry group $G_i$
 such that all the  $G_i-$orbits have the same minimal  dimension  $k. $
When $k=0,$ the  equation $(E_{\a f}^0)$ is critical  and when
$k>0,$ one has  $\frac{n+2-k}{n-2-k}> \frac{n+2}{n-2}$ and
$(E_{\a f}^k)$ turns out to be overcritical.
 The study of equation $(E_{\a f }^k)$ is strongly related to the
 notion of first and second best constants in the Sobolev inequalities
 presented in  section \ref{prelim}.
 The first best constant   appears to be of importance  in
existence results  and the second in  multiplicity results.\\

\section{\bf Preliminaries}\label{prelim}
\def\theequation{2.\arabic{equation}}\makeatother
\setcounter{equation}{0}

Let $(M,g)$ be a compact Riemannian $n$-manifold, $Is(M,g)$ its
isometry group ($Is(M,g)$  is a compact Lie group),
and  $G$   a subgroup of  $Is(M,g).$
By taking its closure $\bar G$ for the standard topology,
 we can assume that $G$ is compact.
 We note for any $p\in [0,+\infty],$
\begin{eqnarray*}
C^p_G      ( M )
&  =  &  \{ u\in C^p(M), \forall \sigma\in G,u\circ \sigma = u   \}\\
H_{1, G}^2 ( M )
 &  =  &  \{ u\in\H,      \forall \sigma\in G,u\circ \sigma = u   \}
\end{eqnarray*}
where the Sobolev space $\H$ is the completion of $C^\infty(M) $ with
respect to the norm
$
\|u\|_{H_1^2}^2 = \|\n u \|_2^2 +\|u\|_2^2.
$
When no confusion is possible, we write $C^p_G, H_1^2, H_{1,G}^2$ instead
of  $C^p_G(M), H_1^2(M), H_{1,G}^2(M).$ If $n-k  >  2,$ 
we let $2^\sh =  \frac{2(n-k)}{n-2-k},$ and
 Hebey-Vaugon \cite{HVsobosym} proved that  for any $1\leq  q \leq 2^\sh,$
  the embedding $H_{1,G}^2 \subset L^q$
is continuous, and compact  if $q<2^\sh.$
For $p< 2^\sh -1,$  compactness of the embedding $H_{1,G}^2 \subset
L^{p+1}$
 implies,  thanks to the variational method,  that there exists a
  $C^\infty_G$ solution for  the equation
\begin{equation}\tag{$E_p$}
\D u + \a u = f u^{p}, \qquad u>0
\end{equation}
where $\D = -\mbox{ div}  {(\n \  )} $ is the $g-$Laplacian, $\a\in
\R^{+ *},$ and $ f\in C_G^\infty$ is positive. When  $p=2^\sh-1,$
the existence of solutions is more difficult to obtain because of
lack of compactness.

 For convenience in what follows, we recall  some  results about the
action of an isometry group $G$ on a compact  manifold. We refer to
Bredon
 \cite{bredon}, Gallot-Hulin-Lafontaine \cite{GHL} and Hebey-Vaugon
 \cite{HVsobosym} for more details.
  Since we can choose $G$ compact, for any $x\in M, $
$O_x^G =\{\sigma(x), \sigma\in G \}$  the $G$-orbit of $x$ is a
compact submanifold of $M$ and $S_x^G=\{\sigma\in G, \sigma(x)=x\}$
 the  isotropy group    of  $x$ is a Lie group of $G.$
A $G$-orbit $O_x^G$ is  principal if for any $y\in M, \ S_y^G$
possesses a subgroup which is conjugate to $S_x^G.$
Principal orbits are of maximum dimension but the converse is false in general.
Let $\Omega$ be the union of all principal orbits. Then $\Omega$ is a
dense open subset of $M,$ and $\Omega/G$ is a quotient manifold. More
precisely, if $\pi $ is the associated submersion, then $(\pi, \Omega,
\Omega/G)$ is a fibration where each fiber is a $G$-orbit.
Note that if all $G$-orbits are principal,
 there exists a unique  manifold structure on the topological
space $M/G   $ and   the metric $g$ induces  a quotient metric $\t
g$ on $M/G$ such that $\pi_G : M\rightarrow M/G$ is a Riemannian
submersion.

We consider here  $\m C^\infty_G $ solutions of $(E_p)$ for
$p=2^\sh-1. $ The equation is written as
\begin{equation}\tag{$E_{\a f}^k$}
\D u + \a u = f u^\frac{n+2-k}{n-2-k},\qquad u>0.
\end{equation}
When $k>0,$ namely when  there is no finite $G$-orbit,
then $\frac{n+2-k}{n-2-k}>\frac{n+2}{n-2}$ and
 $(E_{\a f}^k)$ is, in some sense,  overcritical.
The study of  $(E_{\a f}^k)$ is strongly related to the
problem of the attainability of sharp constants  in functional inequalities
associated with  the continuous embedding  $H_{1,G}^2 \subset L^{2^\sh}.$
Following Faget  \cite{zoe04}, we introduce two assumptions
 $\mHun $ and $\mHdeux $ given by :

\vspace{2mm} \textbf{$(\mathcal{H}_1)$} :  for any orbit $O_{x_0}^G$
of minimum dimension  $k$ and minimum  volume  $ A,$
there exists $H$ a  subgroup of $Is(M,g)$ and  $\d>0$ such that \\
 i)  in  $\m O_{x_0,\d} = \{x\in M / d_g(x,O_{x_0}^G) <\d \}$,
all $H$-orbits are principal,\\
ii) for any  $ x\in \m O_{x_0,\d}, \  O_x^{H} \subset O_x^G$  and
 $O_{x_0}^{H}=O_{x_0}^{G} $,\\
iii) for any  $ x\in \m O_{x_0,\d},
 \  A = vol_g O_{x_0}^G  \leq vol_g O_{x }^H$.
\vspace{2mm}

\noindent and

\vspace{2mm} \textbf{$(\mathcal{H}_2)$} :  for any orbit $O_{x_0}^G$
of minimum dimension  $k$ and minimum  volume  $ A,$ there exists
$H$ a normal
subgroup of $G$ and  $\d>0$ such that \\
 i)  in $ \m O_{x_0,\d} = \{x\in M / d_g(x,O_{x_0}) <\d \}$,
 all  $H$-orbits are principal,\\
ii)  $O_{x_0}^H = O_{x_0}^G$.\\
iii) for any $ x\in \m O_{x_0,\d}, \  x \not\in O_{x_0}^G,
\ dim O_{x}^G >k= dim O_{x_0}^G,$\\
iv)  for any $ x\in \m O_{x_0}$, $x $ is a critical point of the
function $v_H(y)=vol_g O_{y}^H.$

\vspace{2mm}

\noindent  Faget  \cite{zoe04} shows that :


\vspace{2mm}\noindent  {\bf Theorem F }[\textbf{Faget
\cite{zoe04}}]\,\, {\it Let $(M,g)$ be a compact Riemannian
$n$-manifold,  $G$ a compact subgroup of $Is(M,g),$  $k$ the minimum
$G$-orbit dimension, and $A$ the minimum volume of $G$-orbits of
dimension $k.$
 Assume that   $n-k>2.$
If at least one of the assumptions $\mHun$ or
$\mHdeux$ holds true, then there exists   $B>0$ such that for any $
u\in H_{1,G}^2,$
\begin{equation}\label{Isobosanse}
\|u\|_{2^\sh}^2
\leq \frac{K_{n-k}}{A^\frac{2}{n-k}}
     \left[ \|\nabla u\|_2^2 + B \|u\|^2_2 \right],
\end{equation}
where
$K_{n-k} = \frac{4}{(n-k)(n-2-k)\w_{n-k}^{2/(n-k)}}, $ and
$\w_{n-k}$ is the volume of the standard sphere $(S^{n-k}, h_{n-k}).$
The value  $K_{n-k} \  A^{-\frac{2}{n-k}}$ is the best possible
in (\ref{Isobosanse}),  i.e. the smallest constant such that
 (\ref{Isobosanse}) holds true for all $u\in H_{1,G}^2.$} \\

When assumptions $\mHun$ or $\mHdeux $ hold true for a subgroup $H$
we use in the sequel the following notations : $\pi_H$ is the
canonical submersion $\m O_{x_0, \d} \rightarrow \m O_{x_0, \d}/H$
and $\t g$ is the quotient metric induced by $g$ on  $\m O_{x_0,
\d}/H$
 such that  $\pi_H$ is a Riemannian submersion.
For any $x\in \m O_{x_0,\d}, $ we note
 $\t x = \pi_H (O_x^H)$
and $\t v_H$ the function defined for any $y\in \m O_{x_0, \d}/H$ by
$\t v_H (y) = vol_g (\pi_H^{-1}(y)). $

 When inequality (\ref{Isobosanse}) holds true,
  we define the second best constant by
$$
\BoG := \inf \{ B>0,   \forall u \in H_{1,G}^2,\ (\ref{Isobosanse})
 \mbox{ is valid with } B\}.
$$
If (\ref{Isobosanse}) holds true,
we can take $B= \BoG $ in  (\ref{Isobosanse}),  so that for any
$u\in H_{1,G}^2,$
\begin{equation}\tag{$I_{S}^{G,opt}$}
{\|u\|}_{2^\sh}^2
 \leq \frac{K_{n-k}}{A^\frac{2}{n-k}}
      \left( \| \n u \|_2^2 + \BoG\  \| u \|^2_2\right).
\end{equation}
This inequality is  optimal with respect to the first and to the
second constants, i.e. none of them can be improved.
When no confusion is possible we write $B_{0,G}$ instead of $\BoG.$
Note that Hebey-Vaugon \cite{HV} proved earlier that when  $G=\{Id\},$
then $(I_{S}^{Id, opt})$ holds true
on every compact Riemannian $n$-manifold, $n\geq 3.$
As a remark,  $(I_S^{G,opt})$ is  true
  if all $G$-orbits are principal of constant volume,
 since we can take  $H=G$ in  $\mHun$.  We then easily see that
\begin{equation}\label{Bofiniteppalorbits}
B_{0,G}(M,g) = B_{0,Id}(M/G,\t g).
\end{equation}

 Now we discuss the role of the first best constant in
$(I_{S}^{opt})$
 with respect to the existence of  solutions of $(E_{\a f }^k).$
 $G$-invariant  solutions of  $(E_{\a f}^k)$  can be  obtained by
the variational method by minimizing   $I$ on $\mathcal P$ where :
$$
I(u) = \frac{\|\n u\|_2^2 + \a  \|u\|_2^2}
            {\left(\int_M f |u|^{2^\sh}\ dv_g   \right)^{2/2^\sh}},
$$
and
$$
\mathcal P =\left\{ u\in H_{1,G}^2,\  \int_M f |u|^{2^\sh}\ dv_g>0  \right\}.
$$
We note
$
\U_G := \inf_{u\in \mathcal P} I(u).
$
The main difficulty is the lack of compactness coming from the critical
exponent $2^\sh,$  but this is by now a classical problem.
It was firstly solved for the Yamabe problem by working with
subcritical  exponent and then by passing to the limit exponent.
Faget \cite{zoe03} proves that
\begin{equation}\label{majUpsilon}
\U_G \leq \frac{A^\frac{2}{n-k}}{K_{n-k}\ (\max f)^{2/2^\sh}},
\end{equation}
and that, if
\begin{equation}\label{<}
\U_G < \frac{A^\frac{2}{n-k}}{K_{n-k}\ (\max f)^{2/2^\sh}},
\end{equation}
 then there exists a  solution $u\in C^\infty_G$ for $(E_{\a f}^k)$
such that $\U_G = I(u). $ Such a solution is
said to be $G$-minimizing.
Let $(E_\a^k)$ be $(E_{\a f}^k)$ when $f=1.$
Propositions 1 and 2 below follow from the work of Faget \cite{zoe03}.

\begin{proposition}\label{existencesiISGopt}
{\it Let   $(M,g)$ be a compact Riemannian $n$-manifold ,  $n\geq
3,$
      $G$ an isometry group,
      $k$ be the minimum $G$-orbit dimension.
Assume that $n-k>2 $ and that $(I_S^{G, opt})$ holds true. If $
\a \in ]0,  B_{0,G} [, $ then there exists a  $ C^\infty_G$ and
 $G$-minimizing solution for the equation  $(E_{\a}^k).$}
 \end{proposition}

\noindent
\proof
By  the definition of $B_{0,G},$ the strict inequality  $(\ref{<}) $
 holds true, and we can apply the results in Faget \cite{zoe03}.$\hfill\blacksquare$


 \begin{proposition}\label{testueG}
{\it Let   $(M,g)$ be a compact Riemannian $n$-manifold , $n\geq 4,$
      $G$ an isometry group,
      $k$ be the minimum $G$-orbit dimension,
and   $A$ be the minimum volume of $G$-orbits of dimension $k.$
 Assume that   $n-k\geq4.$ Let $x_0\in M$ such that $dim O_{x_0}^G= k $
 and $vol_g O_{x_0}^G= A$ and let $f\in C^2_G$ maximal  at $x_0.$
 Assume that  one of the assumptions $\mHun$ or
$\mHdeux$ holds true for a subgroup $H.$
With the notations introduced above, if
\begin{equation}\label{CSueG}
\left\{
\begin{array}{l}
 (n-4-k)\ \D_{g} f (x_0)= 0 \qquad \\
\\
\a < \frac{n-2-k}{4(n-1-k)}
     \left( \frac{3 \D_{\t g} \t v_H (\t x_0)}{ A }
        +   S_{\t g}(\t x_0) \right),
\end{array}
\right.
\end{equation}
then there exists a  $G$-minimizing $C^\infty_G $   solution for the
equation $(E_{\a f}^k)$.}
\end{proposition}

\noindent\proof For any $\e>0, $ let $\t u_\e$ be  defined on  $\m
O_{x_0, \d}$ by $ \t u_\e =(\e +\t r^2 )^{1-N/2} - (\e +
\d^2)^{1-N/2} $ where $\t r = d_{\t g}(.,\t x_0)$ and $ N=n-k$. We
set  $u_\e = \t u_\e \circ \pi_H,$  and after lengthy computations,
we get that
\begin{equation*}
\begin{array}{ll}
     I (u_\e )
  \leq  \frac{A^{2/N}}{K_{N}  f( x_0)^{2/2^\sh}} \\
 \times  \left[1 + \frac{\e}{N(N-4)}
                \left(    \frac{   \a \ 4(N-1)     }{N-2   }
                       +  \frac{ (N-4) \D_{ g} f ( x_0 )}
                               {2   f( x_0)      }
                       -  \frac{3 \D_{\t g} \t v_H (\t x_0)}{A}
                       -  S_{\t g}(\t x_0)
                \right)  + \circ(\e)
         \right]
&\mbox{if } N>4\\
&\\
\times  \left[  1 + \frac{\e \ln \e }{8}
               \left(    S_{\t g} (\t x_0)
                      +  \frac{3 \D_{\t g} \t v_H (\t x_0) }{A}
                      -  6 \a
                \right)   + \circ(\e \ln \e)
       \right]
&\mbox{if } N=4.
\end{array}
\end{equation*}
Thanks to (\ref{CSueG}), inequality (\ref{<})
holds true  and  we can apply the results in Faget \cite{zoe03}.
Proposition  \ref {testueG} is proved.$\hfill\blacksquare$

\vspace{3mm} Now we briefly discuss estimates on $\BoG.$ At the
moment, the only  compact Riemannian manifold where one knows its
explicit value  is the standard sphere $(S^n,h_n)$ when no isometry
invariance is requiered,
 i.e. when $G=\{Id\}. $
Noting $B_0$ instead of $B_{0, Id},$ one has that
\begin{equation}\label{BoSn}
B_{0}(S^n,h_n) = \frac{n(n-2)}{4}.
\end{equation}
Lower bounds for $\BoG$ have recently been  obtained by Faget \cite{zoe04} :
on a compact Riemannian $n$-manifold, $n\geq 4,$
with the same $G, k, A$ and notations as above, if
$n-k>4$ and if $\mHun$ or $\mHdeux$ holds true, then
\begin{equation}\label{minBoG}
B_{0,G}(M,g)
\geq \max \left\{
     \frac{A^\frac{2}{n-k}}{V_g^\frac{2}{n-k} K_{n-k}},
     \frac{n-2-k}{4(n-k-1)}
\left( S_{\t g}(\t x_0) + \frac{3  \D_{\t g} \t v_H (\t x_0)}{A} \right)
          \right\}
\end{equation}
where $V_g$ is the volume of $(M,g).$
We do not know yet upper bounds for $\BoG$ in the general case.
Hebey-Vaugon \cite{HVmult} computed upper
bounds on specific conformally flat manifolds.
On $(S^1(t)\times S^{n-1}, h_1\times h_{n-1}), $ with $t>0, n\geq 3 $
and when no isometry invariance is requiered, i.e. $G=\{Id\}  :$
\begin{equation}\label{Bogcerclesphere}
\frac{(n-2)^2}{4}
\leq B_{0}(S^1(t) \times S^{n-1}, h_1\times h_{n-1})
\leq \frac1{4 t^2 }  + \frac{(n-2)^2}{4}.
\end{equation}
Note that this approximation is optimal when $t\rightarrow \infty.$
On the quotient manifold $(S^n/G,\t g),\\   n\geq 3,$ where
$G\subset O(n+1)$  is a cyclic group of order $A$ and acts freely on
$S^n$ and $\t g$ is the  quotient metric induced by $h_n, $
\begin{equation}\label{Bogquotientsphere}
\frac{A^{2/n} n (n-2)}{4}
\leq B_{0}(S^n/G,\t g)
\leq \left( 1+\frac{A^2}{4} \right)
     \left( \frac{n+1}{2}  \right) -1 + \frac{n(n-2)}{4}.
\end{equation}
As we will see, these estimates on $B_{0,G},$
especially the upper bounds, are fundamental in the problem
of multiplicity of solutions.\\

\section{\bf Multiplicity results 1}\label{first}
\def\theequation{3.\arabic{equation}}\makeatother
\setcounter{equation}{0}

Assuming that there exists two invariant solutions for $(E_{\a
f}^k), $ we give general   conditions to separat the energies in
Theorems \ref{mult1}.a and \ref{mult1}.b. Then we illustrate these
theorems on specific examples where existence and multiplicity are
compatible. We postpone the proof of  Theorems \ref{mult1}.a and
\ref{mult1}.b to section \ref{proof1}.

\noindent{\bf Theorem 1a}
 \,\,{\it Let $(M,g)$ be a compact Riemannian  $n$-manifold, $n>4,$
 $G_1$ and $G_2$ be  two isometry groups  such that the
minimum dimensions of $G_1$- and $G_2$-orbits are the same. We
denote by  $k\geq 0$ this common minimum orbit dimension, and let
$A_i>0$
 be the minimum volume of $G_i$-orbits of dimension $k,  i\in\{1,2\}.$
 We suppose that $n-k>2,  \   A_1 < A_2$
and that $(I_S^{G_2,opt})$ is valid.
Assume that for  $\a\in \R^*_+$ and $f\in C_{G_1\cup G_2}^\infty$  positive
there exist  two solutions of $(E_{\a f}^k)$ :
 $u_1\in C^\infty_{G_1}$ which is $G_1$-minimizing and
$u_2\in   C^\infty_{G_2}$ which is $G_2$-minimizing.
 If
 \begin{eqnarray}
i)   \qquad   \a   & \leq &   \Bodeux, \label{ia}\\
ii)  \qquad   \a   & \geq &   \frac{n(n-4)}{(n-2)^2}\  \Bog, \mbox{  and }
\label{iia} \\
iii) \qquad   \a   &   >  &
     \Bodeux
 -   \left[  \left(  \frac{A_2}{A_1}  \right)^\frac{2}{n-k}    -1  \right]\
     \frac{A_2^{-\frac{ 4 }{(n-k)(n-2)}}\    K_{n-k}^{\frac{2}{n-2}}}
          {V_g^\frac{2(n-2-k)}{(n-k)(n-2)}\  K_n^\frac{n}{n-2} }\nonumber\\
 &&    \qquad \left( \frac{n(n-4)}{(n-2)^2}  \right)^\frac{n}{n-2}
        \left( \frac{\max f}{  <f>}\right)^\frac{2(n-2-k)}{(n-k)(n-2)},
  \label{iiia}
\end{eqnarray}
where $<f>$ stands for the average value of $f,$
then $ \mE(u_1) < \mE(u_2). $ In particular,  $u_1$ and $u_2$ are
distinct.}\\

 With similar global arguments, and basically only one technical
variation in the proof, we can prove  a slightly different result
:\\

\noindent{\bf Theorem 1b}\,\, {\it Let $(M,g)$ be a compact
Riemannian $n$-manifold, $n>4,$
 $G_1$ and $G_2$ be two  isometry groups  such that the
minimum dimensions of $G_1$- and $G_2$-orbits are the same. We
denote by  $k\geq 0$ this common minimum orbit dimension, and let
$A_i>0$ be the minimum volume of $G_i$-orbits of dimension $k,
i\in\{1,2\}.$ We suppose that $n-k>4,   \  A_1 < A_2$ and that
$(I_S^{G_2,opt})$ is valid. Assume that for  $\a\in \R^*_+$ and
$f\in C_{G_1\cup G_2}^\infty$  positive, there exist  two solutions
of $(E_{\a f}^k)$ :
 $u_1\in C^\infty_{G_1}$ which is $G_1$-minimizing  and
$u_2\in C^\infty_{G_2}$ which is $G_2$-minimizing.
If
 \begin{eqnarray}
i)   \qquad   \a   & \leq &   \Bodeux, \label{ib}\\
ii)  \qquad   \a   & \geq &   \frac{(n-k)(n-4-k)}{(n-2-k)^2} \Bodeux,
\mbox{ and } \label{iib}  \\
iii) \qquad   \a   &   >   &
     \Bodeux
 -   \left[   \left(  \frac{A_2}{A_1}  \right)^\frac{2}{n-k}  -1  \right]
                      \frac{A_2^\frac{2  }{n-k}}
                           {  V_g^{\frac{2}{n-k}} K_{n-k}}\nonumber\\
 &&    \qquad \left(\frac{(n-k)(n-4-k)}{(n-2-k)^2} \right)^\frac{n-k}{n-2-k}
              \left(\frac{\max f}{  <f>}  \right)^\frac{2}{n-k},
\label{iiib}
\end{eqnarray}
where $<f>$ stands for the average value of $f,$
then $ \mE (u_1) < \mE(u_2). $ In particular,   $u_1$ and  $u_2$ are
distinct.}\\

  As a remark, if in Theorems 1.a and 1.b,
  one of the solutions $u_1$ or $u_2$ satisfies (\ref{<}),
then   inequality $iii)$ is not necessarily strict.
We refer to the  proof of Theorems \ref{mult1}.a and  \ref{mult1}.b
for more details on this claim.
As a remark, the
compatibility of conditions $i), ii) $ and
$iii)$  is not automatic. In our examples,
we  choose $f$ such that the right side in $iii)$ is nonpositive so that
$iii)$  is valid. Then multiplicity holds true  when $\a$ belongs to
the interval defined by $i)$ and $ii)$.
 In the following  Corollary of Theorem \ref{mult1}.a, we give  general
conditions in order to separate energies of an infinity of solutions.\\

\begin{corollary}
{\it Let  $(M,g)$ be a compact Riemannian $n$-manifold with $n\geq
3,$ and $(G_i)_{i\in I}$ a family of isometry groups of $Is(M,g)$
such that for any $i\in I,$  $(I_S^{G_i,opt})$ is valid. For any
$i\in I,$ let $ k_i$ be the minimum  dimension of $G_i-$orbits, and
$A_i$ be  the minimum volume of  $G_i$-orbits of dimension $k_i.$ We
assume that $ \forall i\in I,    k_i=k . $
 Given $\a\in \R^{+*},$ and $f\in C^\infty_{\cup_{i\in I}G_i } $  positive,
we suppose that for any  $i\in I$, there exists a
$G_i$-minimizing solution  $u_i\in C^\infty_{G_i}$  for  $(E_{\a f }^k)$.
If $\a\in \left[ \frac{n(n-4)}{(n-2)^2} \Bog;\ \min_{i\in I} (B_{0,G_i})
\right]$ and if
for any $ i\in I$ and $  j\in I$ such that  $A_j < A_i$  we have that
 \begin{equation*}
 \left( \frac{A_i}{A_{j}}  \right)^\frac{2}{n-k}
> 1 + (B_{0,G_i} - \a)
\frac{K_n^\frac{n}{n-2}}{K_{n-k}^\frac{2}{n-2}}A_i^\frac{4}{(n-k)(n-2)}
  \left( \frac{(n-2)^2}{n(n-4)}  \right)^\frac{n}{n-2}
 \left(\frac{\int_M f \ dv_g}{\max  f}\right)^\frac{2(n-2-k)}{(n-k)(n-2)},
\end{equation*} 
then $ \mE(u_j) <  \mE(u_i). $}\\
\end{corollary}

 Now we discuss specific examples. The two  first examples concern
critical  equations and the third example concerns overcritical
equations.

\begin{example}\label{intervalnir}
{\rm Let $(S^n,h_n) $ be the standard sphere of odd  dimension  $n\geq5$
and  $G_1$ and $G_2$ be two finite subgroups of
$O(n+1)$ acting freely on $S^n$ of respective cardinal  $1<A_1<A_2.$
Let $f\in C^\infty_{G_1\cup G_2}$
positive and  maximal  at $x_0\in S^n  $ such that
the  derivatives at $x_0$ are zero up to the  order $n-3,$ and let
$<f>$ be the average value of $f.$
If
\begin{equation*}\label{condmultnirf}
    \left( \frac{\max f}{ <f> }    \right)^{2/n}
\geq 
\end{equation*}
\begin{equation*} \left( B_{0,G_2}(S^n,h_n) -\frac{n^2(n-4)}{4(n-2)} \right)
    \left( \frac{(n-2)^2}{n(n-4)}           \right)^\frac{n}{n-2}
           \frac{4 A_2^\frac{4}{n(n-2)}}{n(n-2)}
    \left[  \left(  \frac{A_2}{A_1}  \right)^\frac{2}{n} -1 \right]^{-1}
\end{equation*}
 then there exist at least two $C^\infty$  solutions of different energies for
 the critical  equation
\begin{equation}\tag{$ E_{\a f}^0 $}
\D  u + \a u = f u^\frac{n+2}{n-2}, \qquad u>0,
\end{equation}
when $\a$ belongs to the interval
 $$
\a \in \left[ \frac{n^2(n-4)}{4(n-2)} ;  \frac{n(n-2)}{4}
\right].
$$
One of these solutions is $G_1$-invariant and the other is
$G_2$-invariant.  }
\end{example}

 As a remark, when $\a = \frac{n(n-2)}{4},$ $(E_{\a f}^0)$ is the
Nirenberg  equation and we recover a result of Hebey-Vaugon
\cite{HVmult}.

\vspace{3mm}\noindent{\it Proof of Example \ref{intervalnir}.}
 \,Since $G_i$  acts freely, $S^n/G_i$ is a
 manifold. with a quotient metric induced by $h_n$ noted $\t g_i.$
As  mentioned in section  \ref{prelim},
since the $G_i$-orbits are principal of constant cardinal,
  $(I_S^{G_i, opt}) $ holds true and with
(\ref{Bofiniteppalorbits}) and  (\ref{Bogquotientsphere}), we have that
\begin{equation}\label{minoBoGiSn}
B_{0,G_i}(S^n,h_n) = B_{0}(S^n/G_i,\t g_i) \geq \frac{n(n-2)}{4}.
\end{equation}
  We claim that for $\a\leq  \frac{n(n-2)}{4}$ there exist two
solutions  $u_i\in C^\infty_{G_i}, \ i=1,2, $ $G_i$-minimizing  for
$(E_{\a f}^0).$
The existence for $\a= \frac{n(n-2)}{4}$ is given by Hebey-Vaugon
\cite{HVmult} since  the derivatives of $f$ are zero up to the
 order $n-3.$
Besides, thanks to  Proposition \ref{testueG},
 there exists  $u_i\in C^\infty_{G_i}$ solution of
$(E_{\a f}^0)$ if
$$
\a < \frac{n-2}{4(n-1)} S_{\t g_i} (S^n/G_i) = \frac{n(n-2)}{4}.
$$
Our claim is proved.
Now according to Theorem \ref{mult1}.a, $u_1$ and $u_2$ are distinct if
 the three assumptions (\ref{ia}), (\ref{iia})
 and (\ref{iiia}) hold true. The first condition  (\ref{ia}) holds true if
$\a\leq\frac{n(n-2)}{4}$, thanks to (\ref{minoBoGiSn}). Condition
(\ref{iia}) is stated  here, since  $B_0(S^n,h_n) =
\frac{n(n-2)}{4},  $ as
$$
\a \geq \frac{n^2(n-4)}{4(n-2)}.
$$
With  this lower bound on $\a,$ in order  to get
(\ref{iiia}), it suffices  that
$$
\frac{n^2(n-4)}{4(n-2)} \geq     B_{0,G_2}(S^n,h_n)
 -   \left[  \left(  \frac{A_2}{A_1}  \right)^\frac{2}{n}    -1  \right]\
     \frac{n(n-2)}{4 A_2^\frac{ 4 }{n(n-2)}}\
     \left(        \frac{n(n-4)}{(n-2)^2}  \right)^\frac{n}{n-2}
              \left(   \frac{\max f}{  <f>\ }         \right)^\frac{2}{n}
$$
with a inequality which is not strict, thanks to the remark following Theorem
\ref{mult1}.b. This is exactly  the assumption made on $f.$
 Thus $u_1$ and $u_2$ exists and are distinct when
$$
\a \in \left[ \frac{n^2(n-4)}{4(n-2)}; \frac{n(n-2)}{4} \right]
$$
and  Example  \ref{intervalnir} is proved.$\hfill\blacksquare$

\vspace{3mm} Now we discuss the following example. Here, we apply
Theorem \ref{mult1}.b and  Theorem \ref{mult1}.a does not provide
the result.

\begin{example}\label{S1Sn-1mult1}
{\rm On    $ (S^1(t)\times S^{n-1}, h_1\times h_{n-1})$ with  $n>4,$ and
$t\geq \sqrt{ \frac{n(n-4)}{4(n-2)^2}},$
let  $G_1 = R_1 \times Id_{S^{n-1}}$ and $G_2 = R_2
 \times Id_{S^{n-1}}$ be two isometry groups,  where $R_1$ and $R_2$ are finite subgroups of
 $SO(2)$ with respective cardinal    $A_1<A_2.$
Let $f\in C_{G_1 \cup G_2}^\infty $  positive and
maximal at $x_0 $ with derivatives at $x_0$ equal to   $ 0 $ up to the order
$n-2$  and such that
\begin{equation}\label{fa}
   \left( \frac{\max f}{ <f> }\right)^{2/n}
\geq  \left( \frac{(n-2)^2}{4} + \frac1{4 t^2}   \right)
          \frac{K_n A_2^\frac{4}{n(n-2)} \left( 2\pi t \w_{n-1}\right)^{2/n} }
               {\left(\frac{A_2}{A_1}        \right)^{2/n} -1 }
   \left( \frac{(n-2)^2}{n(n-4)}              \right)^\frac{n}{n-2}.
   \end{equation}
Then there exist at least two  $C^\infty$ solutions of different energies for
the critical equation $( E_{\a f}^0 )$
when  $\a$ belongs to the interval
\begin{equation*}\label{intmult1S1Sn-1}
\a\in
\left[\frac{n(n-4)}{(n-2)^2} \left( \frac{(n-2)^2}{4} + \frac1{4 t^2}\right)\  ;\
      \frac{(n-2)^2}{4}
\right].
\end{equation*}
One of these solutions is $G_1$-invariant and the other is $G_2$-invariant.}\\
\end{example}

\noindent \textit{\it Proof of Example \ref{S1Sn-1mult1}.} The
$G_i$-orbits are finite and principal and thus
$$
\left( S^1(t) \times S^{n-1} \right) /\left( R_i \times Id_{S^{n-1}}\right)
      = S^1\left(\frac{t}{A_i}\right) \times S^{n-1}
$$
with quotient metric  $h_1\times h_{n-1}.$
As already mentioned in section \ref{prelim},   $(I_S^{G_i, opt}) $ holds true and
with (\ref{Bofiniteppalorbits}) and (\ref{Bogcerclesphere})
\begin{equation}\label{minoBoGSn}
B_{0,G_i } \left( S^1(t) \times S^{n-1}, h_1\times h_{n-1} \right)
\geq  \frac{(n-2)^2}{4}.
\end{equation}
  We claim now that for
\begin{equation}\label{a}
\a\leq  \frac{(n-2)^2}{4}
\end{equation}
 there exist two $\m C^\infty$
solutions  for  $(E_{\a f }^0),$   minimizing  for $G_i, \ i\in{1,2}. $
Since the second derivatives of $f$  at $ x_0$
are zero  and  $S_{h_1\times h_{n-1}}\left(S^1(t/A_i)\times S^{n-1}\right) = (n-1)(n-2), $  the
existence condition (\ref{CSueG}) of Proposition \ref{testueG} is written as
$
\a  <  \frac{(n-2)^2}{4}.
$
If $\a=\frac{(n-2)^2}{4}, $  $(E_{\a f }^0)$ is the equation of the  prescribed scalar
curvature problem and
 it is solved  by Escobar-Schoen \cite{EscSch}
on compact conformally flat manifolds if  $f$ has derivatives at
a maximum point which turn out to be  zero up to  the order $n-2.$
 Thus on $(S^1(t/A_i)\times S^{n-1}, h_1\times h_{n-1})$ there exists $\t u_i$  a minimizing solution of
the equation
$$
\D \t u_i + \frac{(n-2)^2}{4} \t u_i =  \t f  \t u_i^\frac{n+2}{n-2},\qquad \t u_i >0.
$$
If $\pi_i : S^1(t)\times S^{n-1} \rightarrow S^1(t/A_i\times S^{n-1})$ is the canonical submersion,
 then $u_i = \t u_i \circ \pi_i $ is a
$G_i$-minimizing   solution of  $(E_{\a f }^0)$ with $\a= \frac{(n-2)^2}{4}.$
Our claim is proved.
Then one has $\mE(u_1)<\mE(u_2) $ if the three assumptions (\ref{ib}), (\ref{iib}) and (\ref{iiib}) of Theorem \ref{mult1}.b
 hold true.
 (\ref{ib})  is valid by (\ref{minoBoGSn}) if  $\a\leq  \frac{(n-2)^2}{4}.$
 (\ref{iib}) holds true thanks to the upper bound in (\ref{Bogcerclesphere}) if
\begin{equation}\label{b}
\a \geq \frac{n(n-4)}{(n-2)^2} \left( \frac1{4 t^2} + \frac{(n-2)^2}{4} \right).
\end{equation}
By  (\ref{Bogcerclesphere}) and  (\ref{fa}), the right side of (\ref{iiib}) is nonpositive.
 Then (\ref{iiib}) is valid.
Thus  existence and multiplicity are compatible if $\a$ satisfies  (\ref{a}) and  (\ref{b})
which is possible  if
$
t \geq \left(\frac{n(n-4)}{4(n-2)^2}   \right)^{1/2}.
$
Example \ref{S1Sn-1mult1} is proved.$\hfill\blacksquare$

\vspace{3mm} Now we discuss an example where there are non constant
dimensions of  orbits and the minimum dimension is $3.$

\begin{example}\label{S1S2Sn-3}
{\rm On   $(S^1(a) \times S^2(b) \times
S^{n-3},h_1\times h_2\times h_{n-3})$ with $  n\geq 10$ and
\begin{equation}\label{condab}
\frac1{4a} < b^2 < \frac{(n-5)^2}{(n-7)(3n^2-26n+57)},
\end{equation}
we  consider the following  isometry groups:
$$
G_1 = Id_{S^1(a)\times S^2(b)} \times O(n-6) \times O(4) \qquad \mbox{  and   }   \qquad G_2 = O(2) \times O(3) \times Id_{S^{n-3}}.
$$
Let $x_0 = (\o,0_{\R^{n-6}},z_0)$ where $\o\in S^1(a)\times S^2(b) $ and $z_0\in S^3$ and
let  $f\in C^\infty_{G_1\cup G_2}$ be a positive function
  maximal at  $x_0$  such that $\D f ({x_0}) = 0$ and
\begin{equation}\label{condfS1S2Sn-3}
  \frac{\max f}{ <f>}
\geq \left( \left(  4 ab^2       \right)^{2/(n-3)}  -1\right)^{-\frac{n-3}{2}}\
        \left(        \frac{(n-5)^2}{(n-3)(n-7)}  \right)^\frac{(n-3)^2}{2(n-5)} .
\end{equation}
Then there exist at least two $C^\infty$ solutions with different energies for
the over critical equation  $(E_{\a f}^3)$
 when $\a$ belongs to the interval
\begin{equation*}
      \left[\  \frac{(n-3)^2(n-7)}{4(n-5)} ;
                \min \left\{\frac{(n-3)(n-5)}{4},
                        \frac{n-5}{4(n-4)}\left(\frac{2}{b^2} +(n-6)(n-7)\right)   \right\}
    \ \right[.
\end{equation*}
One of these solutions is $G_1$-invariant and the other is
$G_2$-invariant.}
\end{example}

\noindent {\it Proof of Example \ref{S1S2Sn-3}.} The $G_2$-orbits
are  $S^1(a)\times S^2(b) \times \{z\},$ where $z\in S^{n-3},$  and
thus they are principal of constant dimension  $3$ and constant
volume  $ 8\pi^2 a  b^2.$ The quotient metric  on
$\left(S^1(a)\times S^2(b)\times S^{n-3}     \right)/G_2 = S^{n-3}$
is $h_{n-3}.$   According to section  \ref{prelim},  $(I_S^{G_2, opt
})$ holds true and with (\ref{Bofiniteppalorbits}) and (\ref{BoSn})
$$
B_{0,G_2}  =  B_0(S^{n-3}, h_{n-3} )= \frac{(n-3)(n-5)}{4}.
$$
The $G_1$-orbit of   $x=( \o, y, z)\in \R^5 \times \R^{n-6} \times
\R^4$ where $\o \in S^1(a)\times S^2(b), $ and $(y, z)\in S^{n-3}, $ is
 $$
O_x^{G_1} = \{\o\} \times S^{n-7}( \|y\|) \times  S^3( \|z\|).
$$
If  $\|y\|\not= 0$ and $\|z\|\not= 0,
  dim O_x^{G_1}= n-4$ is maximum.
For $x_0 = (\o, 0_{\R^{n-6}}, z_0),$ where
 $\o\in S^1(a)\times S^2(b)$ and $z_0\in S^3,  $  we have
 $$
O_{x_0}^{G_1} = \{\o\}\times \{0_{\R^{n-6}}\} \times S^3
$$
and $dim O_{x_0}^{G_1}=3$ is minimum (thus $O_{x_0}^{G_1}$ is not a
principal orbit) and $ vol O_{x_0}^{G_1}=  2 \pi^2.$
We set  $H= Id_{S^1(a)\times S^2(b)\times \R^{n-6}} \times O(4).  $ $H$
 is a normal subgroup of $G_1,$ and
for any $x = (\o, y,z)$ such that $z\not=0,$
$$
O_x^H = \{\o\}\times \{y\} \times S^3(\|z\|),
$$
where  $\|z\|\in ]0,1].$ The maximum volume for  $H$-orbit is
archieved at  $x_0.$ Moreover the  $H$-orbits are  principal and
 $O_{x_0}^H = O_{x_0}^{G_1}.$ If  $x\not\in O_{x_0}^{G_1},$ then
$
O_x^{G_1} = \{\o\}\times S^{n-7}(\|y\|)\times S^3(\|z\|)
$
with  $\|y\|\not= 0$ and $\|z\|\not = 0$
and   $dim O_x^{G_1} = n-4 >3.$
Finally  assumption   $(\mathcal H_2)$ is true with $H$ and
  $(I_{S}^{G_1,opt})$ is valid.
Now in order to get  $G_i$-invariant and -minimizing solutions of $(E_{\a f}^3),$
we use Proposition \ref{testueG}.
The condition (\ref{CSueG})  for  $G_2$   is
\begin{equation}\label{CSG2}
\a  <    \frac{(n-3)(n-5)}{4}.
\end{equation}
For  $G_1,$ we have
$\D_{\t g} \t v_H (\t x_{0})  \geq 0$
and thus (\ref{CSueG})  holds true if
$$
\a < \frac{n-5}{4(n-4)}   S_{\t g} (\t x_{0}).
$$
Thanks to Proposition \ref{SgMSm} below, this  inequality  holds true if
\begin{equation}\label{CSG1}
\a < \frac{n-5}{4(n-4)} \left( \frac{2}{b^2} + (n-6)(n-7) \right).
\end{equation}
Energies of both solutions obtained under conditions
(\ref{CSG2}) and (\ref{CSG1}) are  different  if the three
multiplicity conditions of  Theorem  \ref{mult1}.b
hold true. The first condition (\ref{ib})
 is $
\a \leq \frac{(n-3)(n-5)}{4}
$
and holds true  if  (\ref{CSG2}) does.  The  second one (\ref{iib}) is stated here as
\begin{equation}\label{iimult1b}
\a \geq \frac{(n-3)^2(n-7)}{4 (n-5)}.
\end{equation}
The last condition  (\ref{iiib}) is stated here as \\

 \noindent$\a>$
$$\frac{(n-3)(n-5)}{4}
     \left[ 1 - \left( (4 ab^2)^{2/(n-3)} -1         \right)
                \left(\frac{(n-3)(n-7)}{(n-5)^2}     \right)^\frac{n-3}{n-5}
                \left(\frac{\max f }{ <f> }  \right)^{2/(n-3)}
     \right].
$$
By (\ref{condfS1S2Sn-3}), the right side of this inequality is nonpositive so that  (\ref{iiib}) holds true.
Finally (\ref{CSG2}), (\ref{CSG1}) and (\ref{iimult1b}) guarantee
existence and multiplicity of two solutions
for  $(E_{\a f}^ 3)$ when
$$
\a\in
\left[\frac{(n-3)^2(n-7)}{4(n-5)};\
      \min \left\{ \frac{(n-3)(n-5)}{4}; \frac{n-5}{4(n-4)} \left(\frac{2}{b^2}+ (n-6)(n-7) \right)\right\}
\right[.
$$
This interval is not empty  thanks to  (\ref{condab}). Example
\ref{S1S2Sn-3} is proved.$\hfill\blacksquare$

\vspace{3mm} Proposition \ref{SgMSm} below was used in the above
proof. \vspace{1mm}\begin{proposition}\label{SgMSm} {\it On a
product manifold $(V^{m}\times S^{n-m}, g\times h_{n-m})$ where
$(V^m,g)$ is a compact Riemannian $m$-manifold, we consider the
isometry  groups
$$
G= Id_V \times O(r_1)\times O(r_2),
\qquad \mbox{ and } \qquad
H=Id_V \times Id_{\R^{r_1}} \times O(r_2)
$$
where
$r_1\geq r_2 $ et $r_1+r_2=n-m+1.$
Let $x_0=(\o_0, 0_{\R^{r_1}}, z_0 ) $  with $\o_0 \in V $ and $z_0\in S^{r_2-1}.$
Then assumption  $\mHdeux$ holds true and with the notations used above,  we have that
$$
S_{\t g} (\t x_0) \geq S_g(\o_0) + r_1(r_1-1).
$$}
\end{proposition}

\vspace{2mm} We postpone the proof of Proposition \ref{SgMSm} to
section \ref{annexe}.

\section{\bf Proofs of Theorems 1.a and 1.b}\label{proof1}
\def\theequation{4.\arabic{equation}}\makeatother
\setcounter{equation}{0}

For convenience, we introduce a general inequality : for  $crit>2$
fixed, $\exists P>0, \exists D>0, \forall u\in H\subset \H,$
\begin{equation}\tag{$I_{PD}$}
 \|u\|^2_{crit}
\leq P
     \left[ \|\n u\|_2^2 + D  \|u\|_2^2               \right]
\end{equation}
where $H\subset H_1^2$ is a functional space such that the inclusion
$H\subset L^{crit}$ is critical in sense of being continuous but not
compact. Theorems 1.a and 1.b are direct corollaries of the
following Theorem \ref{mult1}. In order to get Theorem 1.a from
Theorem \ref{mult1},  it suffices to set $H= H_1^2,\  crit
=\frac{2n}{n-2}, \  P=K_n, $ and  $D=\Bog.$ In this case,
$(I_{PD})$ is the optimal Sobolev inequality $(I_S^{Id,opt})$ which
holds true according to Hebey-Vaugon \cite{HV}  on every compact
Riemanian $n-$manifold, $n\geq 3.$ To get Theorem 1.b from Theorem
\ref{mult1},  it suffices to set $H=H_{1,G_2}^2, \  crit = 2^\sh,\
P=K_{n-k}\ A_2^{-\frac{2}{n-k}},$ and  $D=B_{0,G_2}(M,g).$ In this
case, $(I_{PD})$ is the optimal $G_2$-Sobolev inequality
$(I_S^{G_2,opt})$ which holds true according to Theorem [F]  when we
assume   $\mHun$ or $\mHdeux.$

\begin{theorem}\label{mult1}\,
{\it Let $(M,g)$ be a  compact  Riemannian $n-$manifold, $n\geq 3,$
$G_1$ and $G_2$ be two isometry groups such that the minimum
dimensions of $G_1$- and $G_2$-orbits are the same. We denote by
$k\geq 0$ this common minimum orbit dimension,
 and let
$A_i>0$ be the  minimal volume of $G_i$-orbits of dimension  $k,\ i=1,2.$ We
suppose that
$
n-k>2, \ A_1 < A_2
$
and that  $(I_S^{G_2,opt})$ holds true.
Assume that for  $\a\in \R^{+*} $ and $f\in C_{G_1\cup G_2}^\infty$  positive,
there exist two  solutions of $(E_{\a f}^k)$ :
  $u_1\in C^\infty_{G_1}$ which is   $G_1-$minimizing and $u_2\in C^\infty_{G_2} $ which is
  $G_2-$minimizing. If
 \begin{eqnarray}
i)   \qquad   \a   & \leq &   \Bodeux  \label{imult1}\\
ii)  \qquad   \a   & \geq &   \frac{(4-crit) crit}{4} D \label{iimult1}  \\
iii) \qquad   \a   &   >   &
     \Bodeux
 -   \left[  \left(  \frac{A_2}{A_1}  \right)^\frac{2}{n-k}       -1  \right]
     \frac{A_2^\frac{2 - crit }{n-k}             K_{n-k}^{\frac{crit-2}{2}}}
          { V_g^\frac{(crit-2)(n-2-k)}{2(n-k)}   P^\frac{crit}{2} }\nonumber\\
 &&    \qquad \left(        \frac{(4- crit ) crit}{4}  \right)^\frac{crit}{2}
     \left(   \frac{\max f}{  <f> }         \right)^\frac{(crit -2)(n-2-k)}{2(n-k)}
\label{iiimult1}
\end{eqnarray}
then $ \mE(u_1) < \mE(u_2). $ In particular  $u_1$ and  $u_2$ are
distinct.}
 \end{theorem}

\vspace{3mm}\noindent {\it Proof of Theorem \ref{mult1}.} \,Since
$u_i$ is $G_i$-minimizing,  the strict inequality
$\mE(u_1)<\mE(u_2)$ is equivalent to the strict inegality
$\U_{G_1}<\U_{G_2}. $ According to (\ref{majUpsilon}), it suffices
then to prove that
\begin{equation}\label{minoU2}
\frac{A_1^\frac{2}{n-k}}{K_{n-k} (\max f)^{2/2^\sh}} < \U_{G_2}.
\end{equation}
Note that if $u_1$ satisfies (\ref{<}), then the  equality in
(\ref{minoU2}) is  sufficient to get $\mE(u_1)<\mE(u_2).$ Let us now
search for a lower bound for $\U_{G_2}.$ Since $u_2$ is
$G_2$-minimizing and
 with $(I_S^{G_2,opt}),$ we get that
\begin{eqnarray*}
\frac1{\U_{G_2}}
 & \leq & \frac{  (\max f )^{2/2^\sh} }{\U_{G_2}^\frac{n-k}{2}}
          \ \Kdeux
          \left[ \|\n u_2\|_2^2 + B_{0,G_2} \|u_2\|_2^2 \right].
\end{eqnarray*}
Thus
\begin{equation}\label{majl2debutf}
\frac1{\U_{G_2}}   \leq (\max f )^{2/2^\sh}
                        \ \Kdeux
               \left[ 1 + \frac{ B_{0,G_2} - \a }{\U_{G_2}^\frac{n-k}{2}} \|u_2\|_2^2
               \right  ].
\end{equation}
Since by (\ref{imult1}),  $  B_{0,G_2} -\a \geq 0, $ we  search for an upper bound
for   $\|u_2\|_2^2.$ Multiplying  $(E_{\a f} ^ k)$ by $
u_2^{\frac{4}{crit}-1} $  and integrating over $M$ gives  :
$$
  \| \n u_2^\frac{2}{crit}\|_2^2
=        \frac{4}{crit (4-crit)}
\left(     \int_M f u_2^{2^\sh-2+\frac{4}{crit}} \ dv_g
      - \a \int_M u_2^\frac{4}{crit}     \ dv_g
\right).
$$
Then   by H{\"o}lder's inequality
\begin{eqnarray*}
 \int_M f u_2^{2^\sh-2+\frac{4}{crit}} \ dv_g
& \leq &
  \left(\int_M f u_2^{2^\sh}  \ dv_g  \right)^\frac{2^\sh  -2 + \frac{4}{crit}}{2^\sh}
  \left(\int_M     f          \ dv_g  \right)^\frac{2-\frac{4}{crit}}{2^\sh}
\end{eqnarray*}
and by  $(I_{PD})$
$$
                   \| \n  u_2^\frac{2}{crit}\|_2^2
\geq   \frac1{P}   \|u_2^\frac{2}{crit}\|_{crit}^2
  -     D\         \|u_2^\frac{2}{crit}\|_2^2.
$$
In  particular, we have that
$$
\frac1{P}   \|u_2^\frac{2}{crit}\|_{crit}^2
\leq      \frac{4}{crit (4-crit)}
  \left(\int_M f u_2^{2^\sh}  \ dv_g  \right)^\frac{2^\sh-2 + \frac{4}{crit}}{2^\sh}
  \left(\int_M  f             \ dv_g\right)^\frac{2-\frac{4}{crit}}{2^\sh}
$$
$$
 +     \left(D -   \frac{4 \a }{crit (4-crit)}  \right) \int_M u_2^\frac{4}{crit}\ dv_g.
$$
Now  by  (\ref{iimult1}) and since
 $u_2$ is a $G_2$-minimizing solution we obtain that
\begin{equation*}\label{majnorme2f}
\| u_2  \|_2^2
\leq \left(  \frac{4  P}{(4- crit ) crit}  \right)^\frac{crit }{2}
     \U_{G_2}^\frac{crit-2+n-k}{2}
     \left(\int_M            f\ dv_g           \right)^\frac{ crit-2 }{2^\sh}.
\end{equation*}
Reporting this inequality 
in (\ref{majl2debutf})\\
\noindent $$
\frac1{\U_{G_2}} \leq $$
 $$(\max f )^{2/2^\sh} \Kdeux
   \left[ 1 + \left(B_{0,G_2} - \a \right)
      \left(  \frac{4  P}{(4- crit ) crit}  \right)^\frac{crit }{2}
      \U_{G_2}^\frac{ crit -2}{2}
      \left(\int_M  f\ dv_g                     \right)^\frac{ crit-2 }{2^\sh}
   \right]
$$
and  with the upper bound for $\U_{G_2}$ given by
(\ref{majUpsilon}),
\begin{eqnarray*}
\frac1 {\U_{G_2}}
&\leq &(\max f )^{2/2^\sh}  \  \Kdeux \times  \\
&& \left[  1 +
        ( B_{0,G_2} - \a ) A_2^\frac{crit -2 }{n-k}
\frac{P^\frac{crit}{2} }{K_{n-k}^{\frac{crit-2}{2}}}
          \left(        \frac{4}{(4- crit ) crit}           \right)^\frac{crit}{2}
          \left(   \frac{ \int_M  f \  dv_g  }{\max f}         \right)^\frac{crit -2}{2^\sh}
   \right].
\end{eqnarray*}
Note that if $u_2$ satisfies (\ref{<}),  the above
inequality is strict.
Finally thanks to (\ref{minoU2}) we have
  $\mE(u_1)<\mE(u_2)$ if
$$
 ( \max f )^{2/2^\sh} \frac{K_{n-k}}{A_1^\frac{2}{n-k}}
 >  (\max f )^{2/2^\sh}  \Kdeux \times 
$$
$$
    \left[ 1 +
 ( B_{0,G_2} - \a ) A_2^\frac{crit -2 }{n-k}
\frac{P^\frac{crit}{2} }{K_{n-k}^{\frac{crit-2}{2}}}
          \left(        \frac{4}{(4- crit ) crit}  \right)^\frac{crit}{2}
          \left(   \frac{ \int_M  f\ dv_g }{\max f}
          \right)^\frac{crit -2}{2^\sh}
   \right]
$$
or isolating $\a$ and  introducing  $<f> $ the average value  of $  f$ :
\begin{equation*}\label{condfinalf}
\a > B_{0,G_2}
 -   \left[  \left(  \frac{A_2}{A_1}  \right)^\frac{2}{n-k}      -1  \right]
     A_2^\frac{2 - crit }{n-k}
     \frac{K_{n-k}^{\frac{crit-2}{2}}}{P^\frac{crit}{2} }
     \left(        \frac{(4- crit ) crit}{4}  \right)^\frac{crit}{2}
      \left(   \frac{\max f}{ V_g   <f> }
      \right)^\frac{crit -2}{2^\sh}
\end{equation*}
which is exactly (\ref{iiimult1}). Theorem \ref{mult1} is proved.
Note that the remark following Theorem \ref{mult1}.b  is also proved
since if $u_1$ or $u_2$ satisfies
 (\ref{<}), then the previous inequality is not necessarily strict.$\hfill\blacksquare$

\section{\bf Multiplicity results 2}\label{second}
\def\theequation{5.\arabic{equation}}\makeatother
\setcounter{equation}{0}

We provide another  general result for multiplicity in Theorem
\ref{mult2} below. Then we illustrate the Theorem on specific examples.
 We postpone the proof of Theorem \ref{mult2} to section \ref{proof2}.

\begin{theorem}\label{mult2}
 {\it Let $(M,g)$ be a compact Riemannian  $n$-manifold, $n\geq 3,$
 $G_1$ and $G_2$ be  two isometry groups  such that the
minimum dimensions of $G_1$- and $G_2$-orbits are the same. We
denote by  $k\geq 0$
 this common minimum orbit dimension,
 and let  $A_i>0$ be the minimum volume of $G_i$-orbits of dimension $k, \ i=1,2.$
We suppose that
$n-k>2 $ and $   A_1 < A_2,$
and that   $(I_S^{G_2,opt})$ holds true.
Assume that for  $\a\in \R^*_+$ and $f\in C_{G_1\cup G_2}^\infty$  positive,
there exist  two solutions of $(E_{\a f}^k)$ :
 $u_1\in C^\infty_{G_1}$ which is $G_1$-minimizing, and
$u_2\in C^\infty_{G_2}, $ which is $G_2$-minimizing.
If
 \begin{eqnarray}
i)\qquad   \a   & \leq &      \Bodeux, \mbox{ and } \label{imult2} \\
ii)\qquad   \a   & >    &     \Bodeux
 -  \  \frac{A_2^\frac{2}{n-k} - A_1^\frac{2}{n-k}}
            {K_{n-k}\ V_g^\frac{2}{n-k}}
    \  \frac{\inf f }
            {\max f^\frac{2}{2^\sh} <f>^\frac{2}{n-k}}, 
 \label{iimult2}
\end{eqnarray}
where $<f>$ stands for the average value of $f,$
then $ \mE(u_1) < \mE(u_2). $ In particular,  $u_1$ and $u_2$ are
distinct.}
\end{theorem}

 Here again, if  $ u_1$ satisfies  (\ref{<}),
 then inequality $ii)$ is not necessarily strict.
In the following Corollary to Theorem \ref{mult2},  $f=1$ and we obtain three different solutions for  $(E_{\a}^k).$

\begin{corollary}\label{mult2fcste}
{\it Let $(M,g)$ be a compact Riemannian  $n$-manifold, $n\geq 3,$
 $G_1$ and $G_2$ be two isometry  groups such that the
minimum dimensions of $G_1$- and $G_2$-orbits are the same. We
denote by  $k\geq 0$
 this common minimum orbit dimension,
 and let  $A_i>0$ be the minimum volumes of $G_i$-orbits of dimension $k, \ i=1,2.$
 We suppose that
$n-k>2 $ and $   A_1 < A_2,$
and that  $(I_S^{G_1,opt})$ and $(I_S^{G_2,opt})$ hold true. Then :\\
\textbf{1) } If
\begin{equation}\label{condmult2fcste}
    \Bodeux - \frac{A_2^\frac{2}{n-k}}{K_{n-k}\ V_g^\frac{2}{n-k}}
<   \Boun   - \frac{A_1^\frac{2}{n-k}}{K_{n-k}\ V_g^\frac{2}{n-k}}
\end{equation}
 then there exist two solutions of different energies for the equation
\begin{equation}\tag{$ E_{\a}^k$}
\D u + \a u = u^\frac{n+2-k}{n-2-k}
\end{equation}
when $\a$ belongs to the interval
\begin{equation}\label{intervallefcste}
\a\in
\left[ \Bodeux - \frac{A_2^\frac{2}{n-k} - A_1^\frac{2}{n-k}}{K_{n-k}\ V_g^\frac{2}{n-k}} \ ; \
          \min_{i=1,2} B_{0,G_i}(M,g)
       \right[.
\end{equation}
One of these solutions is non constant and $G_1$-invariant, the other is $G_2$-invariant.\\
\textbf{2) } If moreover
\begin{equation}\label{CStriple}
 \frac{A_2^\frac{2}{n-k}}{K_{n-k} V_g^\frac{2}{n-k}} <  \min_{i=1,2} B_{0,G_i}(M,g)
\end{equation}
then the constant solution $\bar u_\a = \a^\frac{n-2-k}{4}$ of $(E_{\a}^k)$ is different from the two  previous solutions given in
 \textbf{1)} when $\a$ belongs to the interval  \\
\begin{equation}\label{inttriple}
\a \in
 \left[ \max \left\{ \Bodeux - \frac{A_2^\frac{2}{n-k} - A_1^\frac{2}{n-k}}
                                    {K_{n-k}\ V_g^\frac{2}{n-k}},
                      \frac{A_2^\frac{2}{n-k}}{K_{n-k} V_g^\frac{2}{n-k}}
             \right\}
 \ ; \
        \min_{i=1,2} B_{0,G_i}(M,g)
 \right[.
\end{equation}}

\end{corollary}

\noindent {\it Proof of Corollary \ref{mult2fcste}.}\, Part
\textbf{1)} is a  corollary of Theorem \ref{mult2} when $f= 1$
 and where existence of  solutions is given by Proposition \ref{existencesiISGopt}.
We have here $\a<B_{0,G_i}.$
In particular (\ref{<}) holds true and by the remark following Theorem \ref{mult2},
 inequality $ii)$ in Theorem \ref{mult2} is not necessarily strict.
Theorem \ref{mult2} claims that the two  solutions have different energies when $\a$ belongs to the interval in  (\ref{intervallefcste}).
 In particular, with (\ref{minBoG}), we have  that
 $$
\a^\frac{n-k}{2} V_g\geq \frac{A_1}{K_{n-k}^\frac{n-k}{2}}.
$$
But $\a^\frac{n-k}{2} V_g$ is the energy of
 constant solution $ \a^\frac{n-2-k}{4}.$ Since $
\mE(u_1)<A_1 K_{n-k}^{-\frac{n-k}{2}},$
we get that  $\mE(u_1) < \mE( \a^\frac{n-2-k}{4})$ and $u_1$ is not constant. Part \textbf{1)} is proved
and
 $$
\mE(u_1)<\mE(u_2) = \U_2^\frac{n-k}{2} < A_2\ K_{n-k}^{-\frac{n-k}{2}}.
$$
Then $\mE(u_2) < \mE(\a^\frac{n-2-k}{4}) $ if $ \a \geq
\frac{A_2^\frac{2}{n-k}}{K_{n-k}\  V_g^\frac{2}{n-k}}. $ This is
compatible with  (\ref{intervallefcste}), thanks to
 (\ref{CStriple}), and part \textbf{2)} is proved.$\hfill\blacksquare$

\vspace{3mm} Now we discuss specific examples. In the  three
following examples, the manifold is  $S^1(t)\times S^{n-1}$ and  we
fix $f\equiv 1.$ The first example concerns the  critical equation
$(E_{\a}^0)$ and the two other examples concern  the overcritical
equation $(E_{\a}^k)$ with $k=1.$ In the first example, we pass from
the Yamabe multiplicity to an interval of multiplicity.

\begin{example}\label{S1Sn-1k=0mult2}
{\rm On $(S^1(t) \times S^{n-1}, h_1\times h_{n-1}),  n\geq3,$
let $G_1 = R_1 \times Id_{S^{n-1}}$ and \\$G_2 = R_2 \times Id_{S^{n-1}}$ be two isometry groups,
where $R_1$ and $R_2$ are finite subgroups of $SO(2)$ with respectif cardinals
$A_1<A_2.$
If
\begin{equation*}\label{CStS1Sn-1mult2}
t> \max \left\{
\frac{A_2 \w_n}{2\pi \w_{n-1}} \left(\frac{n}{n-2}\right)^{n/2} ; \
 \left(\frac{A_2^2\  (2\pi \w_{n-1})^{2/n}}
              {(A_2^{2/n} - A_1^{2/n})\   n(n-2) \w_n^{2/n}}
   \right)^\frac{n}{2(n-1)}
\right\}
\end{equation*}
then there exist at least three $C^\infty $ solutions of different energies for the critical equation
$(E_\a^0)$
when  $\a$ belongs to the interval
\begin{equation}\label{inttripleS1Sn-1}
\begin{array}{l}
\a\in
\left[\max\left\{ \frac{(n-2)^2}{4} + \frac{A_2^2}{4t^2}
                - \frac{\left(A_2^\frac{2}{n} - A_1^\frac{2}{n}\right) n(n-2) \w_n^{2/n}}
                       {4 (2\pi t \w_{n-1})^{2/n}} \ , \
                  \frac{A_2^{2/n} n(n-2) \w_n^{2/n}}
                       {4 ( 2\pi t \w_{n-1})^{2/n}}
           \right\}\ ; \
      \frac{(n-2)^2}{4}
\right].
\end{array}
\end{equation}
One of these solutions is $G_1$-invariant, the other is
$G_2$-invariant and the third one is the constant solution  $\bar
u_\a = \a^\frac{n-2}{4}.$}
\end{example}

As a remark, when $\a=\frac{(n-2)^2}{4}, (E_\a^0)$ is the Yamabe equation on
$S^1(t)\times S^{n-1}$ and we recover a multiplicity result
of Hebey-Vaugon \cite{HVmult}. \\

 \noindent
{\it Proof of Example \ref{S1Sn-1k=0mult2}.} The actions of the
groups  are  already presented in Example \ref{S1Sn-1mult1}. In
particular, $(I_S^{G_i,opt})$ holds true and with
(\ref{Bofiniteppalorbits}) and  (\ref{Bogcerclesphere}) we have that
\begin{equation}\label{encadrement}
 \frac{(n-2)^2}{4}
\leq
B_{0,G_i }  \left( S^1(t) \times S^{n-1}, h_1\times h_{n-1} \right)
\leq  \frac{A_i^2}{4 t^2} + \frac{(n-2)^2}{4}.
\end{equation}
We claim that there exist two  solutions $u_i\in C^\infty_{G_i}$ for $(E_\a^0)$ if
$$
\a \leq \frac{(n-2)^2}{4}.
$$
 The double existence for $\a<\frac{(n-2)^2}{4}$ is indeed given by (\ref{CSueG}).
For $\a=\frac{(n-2)^2}{4}$ this is given by the
 the resolution of the Yamabe problem on $S^1(t/A_i)\times S^{n-1}$ and with similar arguments
 to the one used in Example \ref{S1Sn-1mult1}.
Now Corollary \ref{mult2fcste} guarantees that $u_1, u_2$ and the constant solution have different energies if
 (\ref{condmult2fcste}) and (\ref{CStriple}) hold true. First by (\ref{encadrement}), (\ref{condmult2fcste}) holds true if
$$
\frac{A_2^2}{4 t^2} + \frac{(n-2)^2}{4} - \frac{A_2^{2/n}}{K_n (2\pi t \w_{n-1} )^{2/n}}
<
\frac{(n-2)^2}{4} -   \frac{A_1^{2/n}}{K_n (2\pi t \w_{n-1} )^{2/n}}
$$
namely if
$$
t> \left(\frac{A_2^2  (2\pi \w_{n-1})^{2/n}}{n(n-2) (A_2^{2/n} - A_1^{2/n} ) \w_n^{2/n}}   \right)^\frac{n}{2(n-1)}.
$$
Since $B_{0,G_i} \geq\frac{(n-2)^2}{4}, $ (\ref{CStriple}) holds true if
$$
\frac{A_2^{2/n} n(n-2) \w_n^{2/n}}{4 (2\pi t \w_{n-1})^{2/n}}
<
\frac{(n-2)^2}{4}
$$
namely  if
$$
t> \frac{A_2 \w_n}{2 \pi \w_{n-1}}   \left(\frac{n}{n-2} \right)^{n/2}.
$$
Under these two conditions on $t,$ Corollary \ref{mult2fcste} gives the  triple multiplicity when $\a $ belongs to the interval
 in   (\ref{inttriple}) which contains the interval in  (\ref{inttripleS1Sn-1}, thanks to (\ref{encadrement}).
 Example \ref{S1Sn-1k=0mult2} is proved.$\hfill\blacksquare$

 \vspace{3mm}The next example involves the  Hopf fibration  and concerns
overcritical equations on  $S^1(t)\times S^3.$

\begin{example}\label{Hopfmult2}
{\rm On $(S^1(t) \times S^3, h_1\times h_3),   $ where $t>1,$
let
$$
G_1 = Id_{S^1(t)}  \times \left\{ (\s,\s)/ \s\in SO(2)  \right\} \qquad \mbox{ and } \qquad
G_2 = O(2)        \times   Id_{S^3}
$$
be two isometry groups.
There exist at least two $C^\infty$ solutions of different energies for the overcritical equation
\begin{equation*}\tag{$ E_{\a}^1 $}
\D u + \a u =  u^5, \qquad u>0
\end{equation*}
when $\a $ belongs to the interval
\begin{equation}\label{intHopf}
\a\in
\left[
\frac{3}{4\ t^{2/3}}; \frac{3}{4}
\right[.
\end{equation}
One of these solutions is $G_1$-invariant and nonconstant, the other
is $G_2$-invariant. Besides if $u_2$ is not the constant solution,
then there exist at least  three different solutions when $\a$
belongs to the interval in (\ref{intHopf}). On the other hand, if
$u_2$ is the constant solution, the interval of multiplicity for
$\a$ extends to $[\frac{3}{4\ t^{2/3}}; 1[.$}
\end{example}

\noindent {\it Proof of Example \ref{Hopfmult2}.} The $G_2$-orbits
are $S^1(t) \times \{\o\},$ where $ \o\in S^3.$ Thus
 they are principal of dimension $ 1 $ and constant volume $2 \pi t$ and
we have that $
         \left( S^1(t) \times S^3      \right) / G_2
=    S^3 $ with quotient metric $h_3.$ As already mentioned,
$(I_S^{G_2, opt}) $ holds true and with (\ref{Bofiniteppalorbits})
and (\ref{BoSn})
$$
B_{0,G_2} (S^1(t) \times S^3, h_1\times h_3)= \frac{3}{4}.
$$
The group $\{(\s,\s), \s\in SO(2)\} $ gives the Hopf fibration $S^3 \rightarrow S^2(1/2)$ with fiber $S^1$
and $h_2$ as quotient metric on $S^2(1/2)$.  The $G_1$-orbits are $\{\rho\} \times S^1$ where $\rho\in S^1(t).$
 Thus they are principal of dimension $1$ and constant volume $2\pi$  and we have  $(S^1(t) \times S^3)/G_1 = S^1(t)\times S^2(1/2)$ with quotient metric $h_1\times h_2.$
Here again  $(I_S^{G_1,opt})$ holds true  and
 $$
B_{0,G_1} (S^1(t) \times S^3, h_1\times h_3 ) = B_0(S^1(t) \times S^2(1/2), h_1\times h_2).
$$
Part \textbf{1)} of Corollary \ref{mult2fcste} gives a multiplicity
interval for $\a$  if (\ref{condmult2fcste}) holds true.
 We easily check that
$$
       B_{0,G_2}     - \frac{A_2^{2/3}}{K_3 V_{h_1\times h_3}^{2/3} } =  0
$$
 and that
$$
        B_{0,G_1}     - \frac{A_1^{2/3}}{K_3 V{h_1\times h_3}^{2/3} }
=       B_{0} (S^1(t) \times S^2(1/2), h_1\times h_2)     - \frac{3}{4 t^{2/3}}.
$$
Thus (\ref{condmult2fcste}) becomes here
$$
B_{0,G_1}  > \frac{3}{4t^{2/3}}.
$$
By (\ref{minBoG}) we know that
$
B_{0,G_1} \geq \max \left\{\frac{3}{4t^{2/3}}, 1   \right\}= 1
$ since $t>1,$
and thus (\ref{condmult2fcste}) holds true.  Part \textbf{1)} of Corollary \ref{mult2fcste} guarantees then a double
 multiplicity when $\a$ belongs to the interval in (\ref{intervallefcste}). We easily see that this interval is here
$$
\a\in
\left[
\frac{3}{4\ t^{2/3}}; \frac{3}{4}
\right[.
$$
In this example, (\ref{CStriple}) does not hold true, so part \textbf{2)} of Corollary \ref{mult2fcste} does not apply.
The constant solution $\bar u_\a = \a^\frac1{4} $ exists for any $\a>0.$ If $u_2\not= \bar u_\a$ then there exist at least  three
solutions of different energies when $\a$ belongs to the interval in  (\ref{intHopf}). Now if $u_2 = \bar u_\a$ then
$u_2$ exists for any $\a>0.$  The solution $u_1$ exists when $\a< 1 \leq B_{0,G_1}$ and its energy verifies
$
\mE(u_1)<A_1 K_{3}^{-\frac{3}{2}}.
$
Thus $u_1$ is not constant if
$$
 \frac{A_1}{K_3^\frac{2}{3}} \leq \mE(\bar u_\a) = \a^\frac{3}{2} V_{h_1\times h_3}
$$
namely if $\a \geq \frac{3}{4 t^{2/3}}.$ The interval of double multiplicity is here $[\frac{3}{4 t^\frac{2}{3}}, 1[.$
Example \ref{Hopfmult2} is proved.$\hfill\blacksquare$

\vspace{3mm} The last example involves infinite non principal
orbits.

\begin{example}\label{S1Sn-1mult2} {\rm On $(S^1(t) \times
S^{n-1}, h_1\times h_{n-1})$ with $n\geq 4$ and   $t> \left(
\frac{n-1}{n-3}  \right)^\frac{n-1}{2},$  let
$$
G_1 = Id_{S^1(t)}  \times O(n-2) \times O(2) \qquad \mbox{ and }
\qquad G_2 = O(2)        \times   Id_{S^{n-1}}
$$
be two isometry groups.
There exist at least two  $C^\infty$ solutions of different energies for the overcritical equation
\begin{equation}\tag{$E_\a^1$}
\D u + \a u = u^\frac{n+1}{n-3}
\end{equation}
when  $\a$ belongs to the interval
$$
\a\in
\left[
 \frac{(n-1)(n-3)}{4 t^\frac{2}{n-1}}  ;
 \frac{(n-3)^2}{4}
\right[.
$$
One of these solutions is $G_1$-invariant and nonconstant, the other
one  is $G_2$-invariant.}
\end{example}

\noindent {\it Proof of Example \ref{S1Sn-1mult2}.}
 The group $G_2$ is the same as in Example \ref{Hopfmult2}.  The $G_2$-orbits  are
 $S^1(t)\times\{\o\}$
where $\o\in S^{n-1},$   of
 dimension $1$ and  constant volume $2\pi t. $ The quotient manifold is  $ (S^{n-1}, h_{n-1})$
and $(I_S^{G_2, opt })$ holds true with
$$
   B_{0,G_2}\left(S^1(t)\times S^{n-1}, h_1\times h_{n-1}\right)
 = \frac{(n-1)(n-3)}{4}.
$$
We easily check that
$$
B_{0,G_2}- \frac{A_2^\frac{2}{n-k}}{K_{n-k} V_{h_1\times h_{n-1}}^\frac{2}{n-k}} = 0.
$$
The $G_1$-orbits are sphere products possibly reduced to a point :\\
 $\forall x =  (\o, y, z)\in S^1(t)\times\R^{n-2}\times\R^{2} \subset S^1(t)\times S^{n-1},$
$$
O_{x}^{G_1} = \{\o\}\times S^{n-3}(\|y\|) \times S^1(\|z\|).
$$
For  $x_0 = (\o,0_{\R^{n-2}}, z_0), $
 where $\o\in S^1(t),$ and $ z_0\in S^1,  $ we have that
$
O_{x_0}^{G_1} = \{\o\} \times \{0_{\R^{n-2}}\} \times S^1.
$
Thus $dim O_{x_0}^{G_1} = 1$ is minimum and $vol O_{x_0}^{G_1} = 2\pi.$
Similar  arguments as in the proof of Example \ref{S1S2Sn-3} show that
 $(\mathcal{H}_2)$  holds true if we choose the normal subgroup $H$ of $G_1$ as
$H= Id_{S^1(t) \times \R^{n-2}} \times O(2). $ Thus
 $(I_S^{G_1, opt})$ holds true.
  Now  assumption   (\ref{condmult2fcste}) of Corollary \ref{mult2fcste} becomes
 $$
 B_{0,G_1} > \frac{(n-1)(n-3)}{4 t^\frac{2}{n-1}}.
 $$
By (\ref{minBoG}) we know that
$$
B_{0,G_1} \geq \max \left\{
\frac{(n-1)(n-3)}{4 t^{2/(n-1)}};\
\frac{(n-3)}{4(n-2)} \left( S_{\t g}(\t x_{0}) + \frac{3 \D_{\t g} \t v_H(\t x_{0})}{A_1} \right)
\right\}.
$$
Since $vol O_{x_0}^H= vol O_{x_0}^{G_1}$ is maximal on $H$-orbits we have
$
\D_{\t g} \t v_H (\t x_{0})\geq 0
$
and according to Proposition \ref{SgMSm},
$
 S_{\t g}(\t x_{0})\geq (n-2)(n-3).
$
In particular
$$
B_{0,G_1} \geq
\max \left\{
\frac{(n-1)(n-3)}{4 t^{2/(n-1)}};\
\frac{(n-3)^2}{4}
\right\}
=
\frac{(n-3)^2}{4}
$$
since
$
t> \left( \frac{n-1}{n-3}  \right)^\frac{n-1}{2}.$ Thus  (\ref{condmult2fcste}) holds true.
Finally part \textbf{1)} of  Corollary \ref{mult2fcste}  guarantees a double multiplicity when  $\a$ belongs to the  interval
in (\ref{intervallefcste}) whose endpoints are
$$
B_{0,G_2} - \frac{\left(A_2^\frac{2}{n-1} - A_1^\frac{2}{n-1}\right) }
               {K_{n-1} V_{h_1\times h_{n-1}}^\frac{2}{n-1}}
= \frac{(n-1)(n-3)}{4 t^\frac{2}{n-1}}
$$
and
$$
\min\{ B_{0,G_1}, B_{0,G_2}\}
\geq
\min \left\{\frac{(n-3)^2}{4},  \frac{(n-1)(n-3)}{4 }\right\} = \frac{(n-3)^2}{4}.
$$
Example \ref{S1Sn-1mult2} is proved.$\hfill\blacksquare$

\section{\bf Proof of Theorem \ref{mult2}}\label{proof2}
\def\theequation{6.\arabic{equation}}\makeatother
\setcounter{equation}{0}

The proofs of Theorems \ref{mult1} and \ref{mult2} are similar but
with an important difference  in the way we find an upper bound for $\|u_2\|_2.$
In order to prove Theorem \ref{mult2} it suffices,
as in the proof of Theorem \ref{mult1},  to prove that
\begin{equation}\label{minoU2bis}
\frac{A_1^\frac{2}{n-k}}{K_{n-k} (\max f)^{2/2^\sh}} < \U_{G_2}.
\end{equation}
 We search for a lower bound for $\U_{G_2}$ and similar
arguments as in proof of Theorem \ref{mult1} lead us to inequality
(\ref{majl2debutf})
\begin{equation}\label{majl2debutfmult2}
\frac1{\U_{G_2}} \leq (\max f )^{2/2^\sh}
                  \ \frac{K_{n-k} }{A_2^\frac{2}{n-k}}
             \left[ 1 + \frac{ B_{0,G_2} - \a }{\U_{G_2}^\frac{n-k}{2}} \|u_2\|_2^2 \right].
\end{equation}
Thanks to (\ref{imult2}), $  B_{0,G_2} - \a\geq 0,$ and  we  search now for an upper bound for $\|u_2\|_2.$
Here is where the proof diverges from the proof of Theorem \ref{mult1}.
We obtain with H{\"o}lder's inequality and since $u_2$ is $G_2$-minimizing  that
\begin{equation*}\label{majnorme2f1}
     \int_M u_2^2 \ dv_g
\leq \frac{\U_{G_2}^\frac{n-2-k}{2}}{\min f }
     \left(\int_M  f     \ dv_g \right)^\frac{2}{n-k}.
\end{equation*}
Reporting this inequality in (\ref{majl2debutfmult2}) and isolating $\U_{G_2}$ gives :
 \begin{equation*}\label{minl2fin2}
\U_{G_2} \geq \frac{A_2^\frac{2}{n-k}}{(\max f )^{2/2^\sh} K_{n-k}}
            - \left( B_{0,G_2} - \a \right) \
               \frac {\left(  \int_M  f \ dv_g \right)^\frac{2}{n-k} }
                     {\min_M f }.
\end{equation*}
Finally  (\ref{minoU2bis}), and thus also the strict inequality $\mE(u_1)<\mE(u_2),$ hold true  if
$$
\frac{A_1^\frac{2}{n-k}}{ ( \max f )^{2/2^\sh} K_{n-k}} <
\frac{A_2^\frac{2}{n-k}}{ (\max f  )^{2/2^\sh}  K_{n-k}}
            - ( B_{0,G_2} - \a )
               \frac {\left( \int_M f \ dv_g \right)^\frac{2}{n-k}}
                     {\min f },
$$
or else
\begin{equation*}\label{condfinalf2}
\a > B_{0,G_2}
 -   \frac{A_2^\frac{2}{n-k}-  A_1^\frac{2}{n-k}}{K_{n-k}}
    \  \frac{\min f }
            {(\max f)^{2/2^\sh} \left(V_g  <f> \right)^\frac{2}{n-k} }.
\end{equation*}
The last inequality is not necessarily strict  when $u_1$ satisfies  (\ref{<}).
 Theorem \ref{mult2} is proved.$\hfill\blacksquare$

\section{\bf Proof of Proposition \ref{SgMSm}}\label{annexe}
\def\theequation{7.\arabic{equation}}\makeatother
\setcounter{equation}{0}

 We start with the following Lemma.

\begin{lemma}\label{Sg>Kg}
{\it Let  $(M,g)$ be  a compact Riemannian $n$-manifold, $n\geq 3,$
of constant sectional curvature $K_g(M), $ and
  $G$ be an isometry group such that all $G$-orbits are principal, and  thus of constant dimension $k.$
 Assume that $k<n.$ Then
\begin{equation}\label{formuleSg>Kg}
  S_{\t g} (y) \geq K_g(M) \ (n-k) (n-k-1),
\end{equation}
for all $y\in M/G,$ where $\t g$ is the quotient metric induced by
$g$  on $M/G.$}
\end{lemma}

 As a remark, if the $G$-orbits are finite, the canonical
submersion $\pi : M\rightarrow M/G$ is a local isometry and
inequality
(\ref{formuleSg>Kg}) is an equality.\\

\noindent { \it Proof of Lemma \ref{Sg>Kg}.}\, On  $(M/G,\t g),$
which has dimension $n-k,$ we have the following relation between
the sectional $K_{\t g} $ and the scalar curvature: $S_{\t g}$
\begin{equation}\label{Scalsect}
 S_{\t g} (y) = \sum_{(i,j)\in [1,n-k]^2, i\not=j} K_{\t g}(\t e_i, \t e_j)
\end{equation}
for all $y\in M/G,$
where  $(\t e_1,...,\t e_{n-k})$ is an orthonormal basis of  $T_{y}(M/G).$
 O'Neil's formula links the  sectional curvatures $K_g$ of $M$ and $K_{\t g }$ of $  (M/G)  $ by
$$
      K_{\t g}(\t e_i, \t e_j)
 =    K_g (e_i, e_j) + \frac{3}{4} \big|[e_i\  e_j]^v\big|^2
\geq  K_g (e_i, e_j)
$$
where $e_i = \left(d\pi_x \setminus_{(Ker d\pi_x)^\perp} \right)^{-1} (\t e_i) \in \left(Ker\ d\pi_x \right)^\perp,$
and where  $[e_i\ e_j]^v \in Ker\ d\pi_x$ is the vertical composant of  $[e_i\ e_j]\in T_x(M). $
Since $K_g$ is constant and with (\ref{Scalsect}), we finally obtain
$$
S_{\t g} (y) \geq K_g(M) \ (n-k)(n-k-1)
$$
and Lemma \ref{Sg>Kg} is proved.$\hfill\blacksquare$

\vspace{3mm} Now we prove Proposition \ref{SgMSm}. \\

\noindent {\it Proof of Proposition \ref{SgMSm}.}\, On the open set
$$
\Omega =  \{x=(\o,y,z)\in V^m\times S^{n-m}, \|z\|\not = 0 \},
$$
 all $H$-orbits are principal and $\mHdeux$ holds true. We have that  $\Omega $ contains
$O_{x_0 }^H = \{\o_0\} \times \{0_{\R^{r_1}}\} \times S^{r_2-1};$
thus there exist an open set $\Omega_1\ni \o_0$ of $V^m$ and an open
set  $\Omega_2\ni \{0_{\R^{r_1}}\} \times S^{r_2-1}$ of $S^{n-m}$
such that
$$
O_{x_0}^H \in \Omega_1 \times \Omega_2 \subset \Omega
$$
and we have
$$
    \left(\Omega_1 \times   \Omega_2 \right)/   H
=   \Omega_1     \ \times \ ( \Omega_2/\ H')
$$
where $H'=Id_{\R^{r_1}}\times O(r_2).$ The metric on  $\left(\Omega_1 \times   \Omega_2 \right)/   H $
 is the  quotient metric  $\t g = g \times \t h_{n-m}$ where $\t h_{n-m} $ is the  quotient metric induced
by $h_{n-m}$ on  $ S^{n-m}/H'.$
Now
$$
\t x_0   = \pi_H \left(\{\o_0\} \times \{0_{\R^{r_1}}\} \times S^{r_2-1}\right)
         = \{\o_0\}\times \{t_0\}
$$
with $t_0 = \pi_{H'} \left(\{0_{r_1}\}\times S^{r_2-1} \right) \in \Omega_2/H'$ and where
$\pi_{H'}: \Omega_2 \rightarrow \Omega_2/H'$ is the canonical submersion. Thus
$$
S_{\t g}(\t x_0)  = S_{g}(\o_0) + S_{\t h_{n-m}}\left(t_0\right).
$$
 Since the  $H'$-orbits are principal on $\Omega_2\subset S^{n-m},$   thanks to lemma
 \ref{Sg>Kg}, and since  $dim \ \Omega_2/H' = n-m-r_2+1= r_1 $ and
$K_{h^{n-m}}(S^{n-m})=1,$ we have
$S_{\t h_{n-m}}(t_0) \geq r_1 (r_1-1). $ Finally
$$
S_{\t g}(\t x_0) \geq S_g(\o_0) + r_1(r_1-1),
$$
Proposition \ref{SgMSm} is proved.$\hfill\blacksquare$

\enddocument